\documentstyle[12pt]{article}
\makeatletter
\def\eqalign#1{\null\vcenter{\def\\{\cr}\openup\jot\m@th
  \ialign{\strut$\displaystyle{##}$\hfil&$\displaystyle{{}##}$\hfil
      \crcr#1\crcr}}\,}
\makeatother

\begin{document}
\def\wt{\widetilde}
\bigskip\bigskip\bigskip
\begin{center}
{\Large\bf Some computable Wiener-Hopf determinants and
polynomials orthogonal on an arc of the unit circle.
}\\
\bigskip\bigskip\bigskip\bigskip
{\large I. V. Krasovsky}\\
\bigskip Technische Universit\"at Berlin\   
Institut f\"ur Mathematik MA 7-2\\
Strasse des 17. Juni 136, D-10623, Berlin, Germany\\
E-mail: ivk@math.tu-berlin.de\\
\bigskip\bigskip
\end{center}
\bigskip\bigskip\bigskip

\noindent
{\bf Abstract.}
Some Wiener--Hopf determinants on $[0,s]$ 
are calculated explicitly for all $s>0$. 
Their symbols are zero on an interval and they are related 
to the determinant with
the sine-kernel appearing in the random matrix theory.

The determinants are calculated by taking limits of Toeplitz
determinants, which in turn are found from the related
systems of polynomials orthogonal on an arc of the unit
circle. 
As is known, the latter polynomials are connected to those orthogonal
on an interval of the real axis. This connection is somewhat extended
here. The determinants we compute originate from the Bernstein-Szeg\H o
(in particular Chebyshev) orthogonal polynomials.

\newpage
\section{Introduction}

Let $\sigma(x)$ be an integrable function on the real line and
\[
K(x)=\delta(x)-{1\over 2\pi}\int_{-\infty}^\infty 
e^{-i\xi x}\sigma(\xi)d\xi.
\]
The Wiener-Hopf operator on $L_2(0,2s)$ corresponding to the symbol 
$\sigma(\xi)$ 
is defined by its kernel $K(x)$
as follows:
\begin{equation}
(W(\sigma)g)(x)=\int_0^{2s}(\delta(x-y)-K(x-y))g(y)dy.
\end{equation}
If $I-W(\sigma)$ is of trace class there exists the determinant (see,
e.g., \cite{BoSi}):
\begin{equation}
\det W=|\delta(x-y)-K(x-y)|_0^{2s}.
\end{equation}

As is well known, Wiener-Hopf determinants $\det W$ are continuous
analogues of the determinants of Toeplitz matrices. 
Let $f(\theta)$ be an integrable function on the unit circle.
Define the  $(n+1)\times(n+1)$ Toeplitz matrix $T_n(f)$ by its
elements as
\begin{equation}
(T_n(f))_{jk}=
{1\over 2\pi}\int_0^{2\pi}e^{-i(j-k)\theta}f(\theta)d\theta,
\qquad j,k=0,1,\dots,n.
\end{equation}
We denote the associated Toeplitz determinant by
\[
D_n(f)=\det T_n(f).
\]

There are many examples of Toeplitz determinants for which an 
explicit expression for all $n$ can be given.
One reason for this is a simple relation between Toeplitz matrices 
and orthogonal polynomials (see below).
The situation in the Wiener-Hopf case is more complicated.
In the present paper we calculate some
Wiener-Hopf determinants whose symbols are zero on an interval.
Note that the known variants of the strong 
Szeg\H o limit theorem, 
which gives large $s$ asymptotics of Wiener-Hopf determinants
(see \cite{BoSi,BW,BWfh}),
are not valid in this case. 
The asymptotics of our Wiener-Hopf determinants have the factor
$e^{-s^2/2}$ in them, and in that, resemble the asymptotics of
Toeplitz determinants with symbols $f(\theta)$ on circular arcs 
found by Widom \cite{WidomArc}. This reflects the fact that, in both
cases, the symbols are zero on intervals (on $(-\alpha,\alpha)$,
$0<\alpha<\pi$, for $f(\theta)$ \cite{WidomArc}, and on $(-1,1)$
for $\sigma(\xi)$ in our examples).    

Similar asymptotics were known for the kernel $K_0(x)=\sin x/\pi x$.
The corresponding Wiener-Hopf determinant 
$\det W_0=|\delta(x-y)-K_0(x-y)|_0^{2s}$ gives the probability
in the Gaussian Unitary Ensemble for an interval of length $2s$ 
(in the bulk scaling limit)
to be free from eigenvalues (see \cite{Mehta}).
The large $s$ asymptotics for it is as follows
\cite{dCM,Dyson,WidomAs,DIZ}:
\[ 
\det W_0=
e^{-s^2/2}s^{-1/4}2^{1/12}e^{3\zeta'(-1)}(1+O(1/s)),
\qquad s\to\infty,
\]
where $\zeta'(x)$ is the derivative of Riemann's zeta function.
Note that appearance of the constant $2^{1/12}e^{3\zeta'(-1)}$ 
has not yet been rigorously justified.

To get our results, we use the following observations:

1) Some Wiener-Hopf determinants can be obtained as the $n\to\infty$ limit
of Toeplitz ones $D_n(f)$ if $f(\theta)$ is allowed to depend
on $n$ \cite{Dyson,BW}. For example, if  $f(\theta)=f_0(\theta)=1$ on the arc 
$\alpha\le\theta\le 2\pi-\alpha$, $\alpha=2s/n$, and otherwise 
$f_0(\theta)=0$, we have $\lim_{n\to\infty}D_n(f_0)=\det W_0$.

2) $D_n(f)$ can be easily obtained given the polynomials
$\Phi_n(e^{i\theta})$ orthogonal on the unit circle with the weight
function $f(\theta)$ (see \cite{Szego} and Lemma 3.1). 

3) Suppose $f(\theta)$ is symmetric ($f(\theta)=f(2\pi-\theta)$).
Then, by the formulas of Szeg\H o, $f(\theta)$ and the polynomials 
$\Phi_n(e^{i\theta})$ are related to a weight function $w(x)$ on the
interval $[-1,1]$ and the orthogonal polynomials associated with it.
A variant of this connection simplified and adopted for
weights $f(\theta)$ on an arc and symmetric $w(x)$ ($w(x)=w(-x)$) 
is Lemma 2.1 below (cf. \cite{Zhedanov}).
The Wiener-Hopf determinant for the kernel $\sin x/\pi x$  corresponds to
\[
w(x)=(1-\gamma^2x^2)^{-1/2},\qquad \gamma=\cos {s\over n}.
\]
Unfortunately, the associated orthogonal polynomials on the interval
$[-1,1]$ are not explicitly known.

On the other hand, many
systems of orthogonal polynomials on an interval have been
extensively investigated. Thus, taking various $w(x)$ for which the
associated polynomials are known, we can hope to
obtain the corresponding Wiener-Hopf determinants. However, for most
known weight functions $w(x)$ the limit $\lim_{n\to\infty}D_n(f)$
does not have proper convergence. So it is only in very
special cases that we get Wiener-Hopf determinants. Those we
calculate correspond to
\begin{equation}
w(x)=\frac{\sqrt{1-x^2}}{1-\gamma^{2r^2}x^2},\qquad r\ge 0.\label{wx}
\end{equation}
In this case the polynomials can be explicitly written in terms of
Bernstein-Szeg\H o polynomials. The corresponding function $f(\theta)$
is asymptotically ($n\to\infty$) $1$ on the unit circle except for a
small (of order $1/n$) neighbourhood of $\theta=0$ (cf $f_0(\theta)$).
Our main result is 

\noindent
{\bf Theorem 1.1} (Bernstein-Szeg\H o  Fredholm determinant.) {\it
Let the kernel
\begin{equation}
K_{BS}(r,z)={\sin z \over \pi z}-{1\over\pi}\int_0^\infty\cos(z\cosh
t)\left({\sinh t \cosh t \over \sinh^2 t +r^2}-1\right)\sinh t dt,
\qquad r\ge 0.\label{KBS}
\end{equation}
Then
\begin{equation}
\left|\delta(x-y)-K_{BS}(x-y)\right|_0^{2s}=e^{-s^2/2-2rs}(\cosh
s + r\sinh s).
\end{equation}

For $r\neq 0$ the kernel can be rewritten as
\begin{equation}
K_{BS}(r,z)={1\over\pi z}\int_0^\infty\sin(z\cosh
t)\frac{r^2\cosh^2 t + (r^2-1)\sinh^2 t}{(\sinh^2 t +r^2)^2}dt,\qquad r>0,
\label{K2}
\end{equation}
and, furthermore, for $r=1$, as 
\begin{equation}
K_{BS}(1,z)=
\int_{|z|}^\infty{J_1(t)\over 2t}dt,\label{Kb}
\end{equation}
where $J_1(t)$ is the 1st order Bessel function.
}

\noindent
{\bf Remark.} As is easy to verify, the symbol
\[
\sigma(r,\xi)=1-\int_{-\infty}^\infty K_{BS}(r,z)e^{i\xi z}dz=
\cases{
0,& $|\xi|<1$,\cr
{|\xi|\sqrt{\xi^2-1}\over \xi^2-1+r^2},& $|\xi|>1$.}
\]

Should we wish to generalize this theorem, we could consider the
weights obtained by multiplying (\ref{wx}) with factors of the form
$(1-x^2+ax^2/n^2)/(1-x^2+bx^2/n^2)$. The corresponding polynomials can
be reconstructed from the Bernstein-Szeg\H o polynomials using the
Christoffel formula (\cite{Szego}, Theorem 2.5).

We also consider the polynomials orthogonal on the arc  
$\alpha\le\theta\le 2\pi-\alpha$ of the unit circle with the weights
$f_1(\theta)=\sin\theta/2$, $f_2(\theta)=\sin^{-1}\theta/2$ (Lemma 2.2).
They are related to the Legendre polynomials.
Note that the polynomials associated with the weight $\sin\theta/2$
were discussed in \cite{GN} where their reflection
coefficients were found.
For a review of polynomials orthogonal on an arc of the unit 
circle, see \cite{GAkh,GN,BL,RL} and references therein.
The determinants $D_n(f_1)$ and $D_n(f_2)$ are given by Lemma 4.1.

\section{Polynomials on an arc}
Let $w(x)$ be a symmetric ($w(x)=w(-x)$) weight function 
on the interval $[-1,1]$ and $0<\gamma<1$. 
Let $P_n(x)=x^n+\dots$, $Q_n(x)=x^n+\dots$, $n=0,1,\dots$ be the systems
of monic orthogonal polynomials associated with the weights $w(x)$
and $w(x)(1-\gamma^2x^2)$, respectively. We have
\begin{equation}
\eqalign{
\int_{-1}^1 P_n(x)P_m(x)w(x)dx=h_n\delta_{mn};\\
\int_{-1}^1 Q_n(x)Q_m(x)w(x)(1-\gamma^2 x^2)dx=h'_n\delta_{mn};\qquad
m,n=0,1,\dots.}
\end{equation}
Since we assumed $w(x)$ to be symmetric the polynomials $P_n(x)$ and
$Q_n(x)$ contain only even (odd) powers of $x$ if $n$ is even (odd).

Furthermore, let $f(\theta)$ be a weight function on the 
unit circle and $\Phi_n(z)=z^n+\dots$, $z=e^{i\theta}$, the
corresponding monic orthogonal polynomials:
\begin{equation}
{1\over 2\pi}\int_{0}^{2\pi} \Phi_n(z)\overline{\Phi_m(z)}f(\theta)d\theta=
{\delta_{mn}\over \chi_n^2};\qquad m,n=0,1,\dots.\label{phi}
\end{equation}

\noindent
{\bf Lemma 2.1} (Connection between polynomials on an arc 
and on an interval.) {\it Let 
$f(\theta)=w(\gamma^{-1}\cos{\theta/2})\sin{\theta/2}$,
if $\alpha\le\theta\le 2\pi-\alpha$, and $f(\theta)=0$ otherwise.
Here $\gamma=\cos{\alpha/2}$, $0\le\alpha\le \pi$. Set
$x=\gamma^{-1}\cos{\theta/2}=(2\gamma)^{-1}(z^{1/2}+z^{-1/2})$, 
$z=e^{i\theta}$.
Then the corresponding polynomials $\Phi_n(z)$ are 
related to $P_n(x)$, $Q_n(x)$ by the following expressions:
\begin{equation}
\eqalign{
\Phi_n(z)={(2\gamma)^{n+1}z^{n/2}\over z-1}
\left(z^{1/2}P_{n+1}(x)-
{P_{n+1}(\gamma^{-1})\over P_{n}(\gamma^{-1})}P_{n}(x)\right);\\
\chi_n^{-2}=(2\gamma)^{2n}{\gamma^2\over\pi}
{P_{n+1}(\gamma^{-1})\over P_{n}(\gamma^{-1})}h_n;\qquad
n=0,1,\dots.}\label{phiP}
\end{equation}

\begin{equation}
\eqalign{
\Phi_n(z)=(2\gamma)^n z^{n/2}\left( Q_n(x)-z^{-1/2}{t_n\over
  t_{n-1}}Q_{n-1}(x)\right);\\
\chi_n^{-2}=(2\gamma)^{2n}{t_n\over\pi t_{n-1}}h'_{n-1};
\qquad n=1,2,\dots;\\
t_k={1\over
  2\pi}\int_{\alpha}^{2\pi-\alpha}z^{k/2}Q_k(x)f(\theta)d\theta;
\qquad k=0,1,\dots.}\label{phiQ}
\end{equation}

\begin{equation}
\chi_0^{-2}=t_0.\label{t0}
\end{equation}
}

\noindent
{\it Remark.} The expressions (\ref{phiP}) were obtained in \cite{Zhedanov}.

\noindent
{\it Proof.} Note first that, since $f(\theta)=f(2\pi-\theta)$, 
the coefficients of $\Phi_n(z)$ are real. This follows from the
definition of the polynomials $\Phi_n(z)$ as a determinant (e.g.,
\cite{Szego}, p.286).
As in the proof of Theorem 11.5 in \cite{Szego}
we establish the formulas for $n=0,1,\dots$:
\begin{eqnarray}
P_n(x)=\frac{z^{-n/2}\Phi_n(z)+z^{n/2}\Phi_n(z^{-1})}{(2\gamma)^n(1-a_{n-1})}=
\frac{\Phi_n(z)+\Phi_n^*(z)}{(2\gamma)^n(1-a_{n-1})z^{n/2}},\label{Pphi}\\
h_n={2\pi \over (2\gamma)^{2n}(1-a_{n-1})\chi_n^2\gamma},\label{hn}
\end{eqnarray}
where $a_n=-\Phi_{n+1}(0)$
and $\Phi_n^*(z)=z^n\Phi_n(z^{-1})$. Note that orthogonality of $P_n(x)$
and $P_m(x)$ with $m$, $n$ of {\it different} parity is obvious from
the symmetry property of the weight. 
As is well known (e.g., \cite{Geronimus}, p.132), the polynomials $\Phi_n(z)$
satisfy the following recurrence relations:
\[
\Phi_{n+1}(z)=z\Phi_n(z)-a_n\Phi_n^*(z),\qquad
\Phi^*_{n+1}(z)=\Phi^*_n(z)-a_n z\Phi_n(z).
\]
These relations substituted into (\ref{Pphi}) where $n$ is replaced by
$n+1$ allow us to express
$P_{n+1}(x)$ as a combination of $\Phi_n(z)$ and $\Phi_n^*(z)$. On the
other hand, (\ref{Pphi}) also gives $P_n(x)$ in terms of 
$\Phi_n(z)$ and $\Phi_n^*(z)$. Eliminating $\Phi_n^*(z)$ from these
expressions gives:
$$
(z-1)\Phi_n(z)= (2\gamma)^{n+1} z^{n/2}\left(
  z^{1/2}P_{n+1}(x)-{1-a_{n-1}\over 2\gamma}P_{n}(x)\right).
$$
Setting here $z=1$, we obtain
\begin{equation}
1-a_{n-1}=2\gamma{P_{n+1}(\gamma^{-1})\over P_n(\gamma^{-1})}.
\end{equation}
The last two equations and (\ref{hn}) yield (\ref{phiP}).

The analogue of (\ref{Pphi}), (\ref{hn}) for $Q_{n-1}(x)$, 
$n=1,2,\dots$ reads:
\begin{eqnarray}
Q_{n-1}(x)=\frac{\Phi_n(z)-\Phi_n^*(z)}
{(2\gamma)^{n-1}(1+a_{n-1})z^{n/2}(z^{1/2}-z^{-1/2})},\\
h'_{n-1}={\pi \over (2\gamma)^{2n}(1+a_{n-1})\chi_n^2}.
\end{eqnarray}
An argument as before with the recurrence relations gives:
$$
\Phi_n(z)= (2\gamma)^n z^{n/2}\left(
  Q_n(x)-{1+a_{n-1}\over 2\gamma}z^{-1/2}Q_{n-1}(x)\right).
$$
Multiplying this equation with $f(\theta)$ and integrating over
$0\le\theta\le 2\pi$, we get by orthogonality of $\Phi_n(z)$ and 
$\Phi_0(z)=1$:
\begin{equation}
1+a_{n-1}=2\gamma {t_n\over t_{n-1}},\qquad
t_n={1\over
  2\pi}\int_{\alpha}^{2\pi-\alpha}z^{n/2}Q_n(x)f(\theta)d\theta.
\end{equation}
Using this, we immediately obtain (\ref{phiQ}). The equation
(\ref{t0}) is obvious from (\ref{phi}).$\Box$

As an example, we shall now present two functions $f(\theta)$ 
for which the corresponding orthogonal polynomials on an arc 
are explicitely given (in terms of Legendre polynomials).

\noindent
{\bf Lemma 2.2} {\it 1) Let $f(\theta)=\sin\theta/2$
if $\alpha\le\theta\le 2\pi-\alpha$, and $f(\theta)=0$ otherwise. Then
$\Phi_n(z)$ are given by (\ref{phiP}) where
\[
P_n(x)=L_n(x);\qquad h_n={2^{2n+1}\over(2n+1){2n\choose n}^2}.
\]
Here
\begin{equation}
L_n(x)={2^n\over{2n\choose n}}\sum_{k=0}^n{n\choose k}{n+k\choose k}
\left(\frac{x-1}{2}\right)^k=x^n+\dots \label{L}
\end{equation}
are the monic Legendre polynomials.

2) Let $f(\theta)=1/\sin\theta/2$
if $\alpha\le\theta\le 2\pi-\alpha$, $\alpha>0$, and $f(\theta)=0$
otherwise. 
Then $\Phi_n(z)$ are given by (\ref{phiQ}), (\ref{t0}), where
\[
Q_n(x)=L_n(x);\qquad h'_n={2^{2n+1}\over(2n+1){2n\choose n}^2}.
\]
}

\noindent
{\it Proof.} Applying Lemma 2.1, we obtain in the case 1 $w(x)=1$
and in the case 2 $w(x)(1-\gamma^2x^2)=1$. Now it only remains
to note that the Legendre polynomials defined by (\ref{L}) satisfy the
orthogonality relation:
\[
\int_{-1}^1
L_n(x)L_m(x)dx=\frac{2^{2n+1}}{(2n+1){2n\choose
    n}^2}\delta_{mn},\qquad n,m=0,1,\dots
\]
$\Box$

\section{Wiener--Hopf determinants}
Consider again monic orthogonal polynomials satisfying relation
(\ref{phi}). A simple but very important fact about the Toeplitz
determinant $D_n(f)$ is that (see, e.g., \cite{Szego}, p.286)
\begin{equation}
D_n(f)=\prod_{j=0}^n \chi_j^{-2}.\label{Dphi}
\end{equation}
Using this we immediately obtain from Lemma 2.1 the following statement:

\noindent
{\bf Lemma 3.1} {\it Let the quantities $f(\theta)$, $\gamma$, $P_k(x)$, $Q_k(x)$, 
$h_k$, $h'_k$, $t_k$ be the same as in Lemma 2.1. Then
\begin{eqnarray}
D_n(f)=2^{n(n+1)}\frac{\gamma^{n^2+3n+2}}
{\pi^{n+1}}P_{n+1}(1/\gamma)
\prod_{j=0}^n h_j;\label{D}\\
D_n(f)=(2\gamma)^{n(n+1)}\frac{t_n}{\pi^{n}}\prod_{j=0}^{n-1} h'_j.\label{D2}
\end{eqnarray}
}

Now we are ready to prove Theorem 1.1. We begin with a particular case
of $r=0$:

\noindent
{\bf Theorem 3.2} (Chebyshev Fredholm determinant.) {\it
Let the kernel
\begin{equation}
K_C(z)={\sin z \over \pi z}-{1\over\pi}\int_0^\infty\cos(z\cosh
t)e^{-t}dt.\label{Kc}
\end{equation}
Then
\begin{equation}
\left|\delta(x-y)-K_C(x-y)\right|_0^{2s}=e^{-s^2/2}\cosh s
\end{equation}

\noindent
Proof.} 
Consider the weight function
\begin{equation}
w(x)={1\over\sqrt{1-x^2}},\qquad x\in[-1,1].
\end{equation}
The corresponding monic orthogonal polynomials $P_k(x)$ are 
the Chebyshev polynomials of the first kind: 
\begin{eqnarray}
P_0(x)=1,\quad h_0=\pi,\\
P_k(x)={1\over 2^{k-1}}\cos(k\;{\rm arccos}\;x),\quad 
h_k={\pi\over 2^{2k-1}}\quad k=1,2,\dots
\end{eqnarray}
Lemma 2.1 gives us now the polynomials $\Phi_k(z)$ orthogonal on the arc 
$\alpha\le\theta\le 2\pi-\alpha$ with the weight function
\begin{equation}
f(\theta)=w\left(\gamma^{-1}\cos{\theta\over 2}\right)\sin{\theta\over 2}=
\frac{\sin{\theta\over 2}}{\sqrt{1-\gamma^{-2}\cos^2{\theta\over 2}}}
.\label{f1}
\end{equation}

Consider $D_n(f)$ for $\alpha=2s/n$, $s>0$ and large $n$. In this case 
$\gamma=\cos{\alpha/2}=1-s^2/(2n^2)+O(n^{-4})$ and (\ref{f1}) takes the
form
\begin{equation}
f(\theta)=\left(1-{s^2\over
    n^2}{\rm cotan}^2\;{\theta\over 2}
\{1+O(n^{-2})\}\right)^{-1/2}.\label{f11}
\end{equation}
As ${\rm arccos}\;\gamma^{-1}=is/n+O(n^{-3})$, we also have
$P_n(\gamma^{-1})=(\cosh s +O(n^{-2}))/2^{n-1}$.

The formula (\ref{D}) gives
\begin{equation}
\eqalign{
D_n(f)=2^{n(n+1)}\left(\frac{(1-s^2/(2n^2)+O(n^{-4}))^{n^2+3n+2}}
{\pi^{n+1}}\right)\left(\frac{\cosh s +O(n^{-2})}{2^n}\right)\times\\
\pi\prod_{j=1}^n\frac{\pi}{2^{2j-1}}=e^{-s^2/2}\cosh s (1+o(1)).}\label{Dc}
\end{equation}
Thus $D_n(f)$ tends to a finite limit as $n\to\infty$. This
indicates existence of the Fredholm determinant
$\lim_{n\to\infty}D_n(f)$. Let us therefore consider the behaviour of
the matrix elements of $T_n(f)$ as $n\to\infty$, that is of
\begin{equation}
I_k={1\over
  2\pi}\int_\alpha^{2\pi-\alpha}e^{-ik\theta}f(\theta)d\theta=
{1\over
  \pi}\int_\alpha^\pi \cos{k\theta}f(\theta)d\theta,\label{Ik}
\end{equation}
where $f(\theta)$ is given by (\ref{f11}). We now split
the last integral into a sum of two: 
one $I_{k1}$, over a small neighbourhood 
$(\alpha,\alpha+\epsilon)$ of
$\alpha$, the other $I_{k2}$, over the interval $(\alpha+\epsilon,\pi)$. 
In the first integral we change the variables $\theta=2xs/n$:
\[
\eqalign{
I_{k1}=\frac{2s}{\pi n}\int_1^{1+\epsilon n/(2s)}\frac{\cos(2sxk/n)}
{\sqrt{1-x^{-2}(1+O(x^2/n^2))}}dx=\\
\frac{2s}{\pi n}\int_1^{1+\epsilon n/(2s)}\cos(2sxk/n)\left(\frac{1}
{\sqrt{1-x^{-2}}}-1\right)dx+\\
 \frac{2s}{\pi n}\int_1^{1+\epsilon n/(2s)}\cos(2sxk/n)dx+O(\epsilon^3).}
\]
Here we added and subtracted cos under the sign of the integral. It
ensures that we can replace the upper integration limit in the
first integral of the last sum by infinity introducing by doing so
the error of order $1/(n^2\epsilon)$. 
After that we change the variable $x=\cosh t$
in the first integral and integrate the second one to get:
\[
I_{k1}=
\frac{2s}{\pi n}\int_0^\infty\cos(2s{k\over n}\cosh t)e^{-t}dt+
{1\over\pi k}\sin(2s{k\over n}+\epsilon k)-
{1\over\pi k}\sin(2s{k\over n})+O\left({1\over n^2\epsilon},\epsilon^3\right),
\]
where for $k=0$ expressions $\sin ak/k$ should be replaced by $a$.
As for $I_{k2}$, we have
\[
I_{k2}={1\over\pi}\int_{\alpha+\epsilon}^\pi\cos(k\theta)d\theta+
O\left({1\over n^2\epsilon^2}\right)=
\cases{
1-{1\over\pi}(\alpha+\epsilon)+O(\{n\epsilon\}^{-2}),& $k=0$,\cr
-{1\over\pi k}\sin(2s{k\over n}+\epsilon k)+O(\{n\epsilon\}^{-2}),
& $k\neq 0$.}
\]

Taking, e.g., $\epsilon=1/n^{2/5}$, we finally obtain
\begin{equation}
I_k=I_{k1}+I_{k2}=
\cases{
1+o(n^{-1}),& $k=0$,\cr
-{1\over\pi k}\sin(2s{k\over n})+
\frac{2s}{\pi n}\int_0^\infty\cos(2s{k\over n}\cosh
t)e^{-t}dt+o(n^{-1}),& $k\neq 0$,}
\end{equation}
uniformly in $k$.
Thus, as a simple analysis shows, 
$\lim_{n\to\infty}D_n(f)=|\delta(x-y)-K_C(x-y)|_0^{2s}$, where 
$K_C(z)$ is given by (\ref{Kc}). In view of equation (\ref{Dc}), we
completed the proof.
$\Box$
  
\bigskip
\noindent
{\it Proof of Theorem 1.1.} 
Take the weight function
\begin{equation}
w(x)=\frac{\sqrt{1-x^2}}{1-\gamma^{2r^2}x^2},\qquad x\in[-1,1],
\qquad r\ge0.\label{w2}
\end{equation}
The corresponding orthogonal polynomials are a particular case of the
Bernstein-Szeg\H o polynomials (\cite{Szego}, Theorem 2.6).
They are given by the expressions:
\begin{equation}
\eqalign{
\wt P_0=\sqrt{2(1-a)\over\pi},\qquad 
\wt P_k(x)=\sqrt{{2\over\pi}}\left((1-2ax^2)\frac{\sin(k+1)\psi}{\sin\psi}+
2ax\cos(k+1)\psi\right),\\
x=\cos\psi,\qquad k=1,2,\dots,}\label{BS}
\end{equation}
where $a=(1-\sqrt{1-\gamma^{2r^2}})/2$.
The polynomials $\wt P_k(x)$ are {\it orthonormal}, that is
$\int_{-1}^1 \wt P_k(x)\wt P_m(x)=\delta_{km}$.
Therefore the monic polynomials and constants $h_k$ are given by the 
formulas:
\[
P_k(x)={\wt P_k(x)\over \kappa_k},\qquad h_k=\kappa_k^{-2},\qquad 
k=0,1,\dots,
\]
where $\kappa_k$ are the coefficients of the highest degree of $\wt
P_k(x)$. Collecting the coefficients of $x^k$ in (\ref{BS}), we obtain:
\[
\kappa_0=\sqrt{2(1-a)\over\pi},\qquad
\kappa_k=\sqrt{{2\over\pi}}2^k(1-a),\qquad k=1,2,\dots.
\]
The weight function of the corresponding polynomials $\Phi_k(z)$ on
the arc $\alpha\le\theta\le 2\pi-\alpha$ 
\begin{equation}
f(\theta)=w\left(\gamma^{-1}\cos{\theta\over 2}\right)\sin{\theta\over 2}=
\frac{\sqrt{1-\gamma^{-2}\cos^2{\theta\over 2}}}
{1-\gamma^{2(r^2-1)}\cos^2{\theta\over 2}}\sin{\theta\over 2}
.\label{f2}
\end{equation}
Note that the functions for $r=0$ and $r=1$ are inverses of each other.
Just as in the previous proof, we
consider $D_n(f)$ for $\alpha=2s/n$, $s>0$ and large $n$. 
Since $a=1/2-rs/(2n)+O(n^{-3})$ we obtain
\begin{equation}
\eqalign{
h_0=\kappa_0^{-2}=\pi\left(1-{rs\over n}+O(n^{-2})\right);\quad
h_k=\kappa_k^{-2}={\pi\over 2^{2k-1}}\left(1-{2rs\over
    n}+O(n^{-2})\right),\\ k=1,2,\dots,\qquad
P_{n+1}(1/\gamma)=2^{-n}(\cosh s+ r\sinh s +o(1)).}
\end{equation}
Substituting this into (\ref{D}) gives
\begin{equation}
D_n(f)=e^{-s^2/2-2rs}(\cosh s+ r\sinh s+o(1)).\label{Db}
\end{equation}
For the elements of the Toeplitz matrix $T_n(f)$ we now have:
\begin{equation}
I_k={1\over
  \pi}\int_\alpha^\pi \cos{k\theta}\frac{\left(1-{s^2\over
    n^2}{\rm cotan}^2\;{\theta\over 2}
\{1+O(n^{-2})\}\right)^{1/2}}
{1+{(r^2-1)s^2\over
    n^2}{\rm cotan}^2\;{\theta\over 2}
\{1+O(n^{-2})\}}
d\theta.\label{2I}
\end{equation}
This integral can be estimated in the same way as the one in the
previous proof. We then obtain the result of the theorem with
\begin{equation}
K_{BS}(r,z)={\sin z \over \pi z}-{1\over\pi}\int_0^\infty\cos(z\cosh
t)\left({\sinh t \cosh t \over \sinh^2 t +r^2}-1\right)\sinh t dt,
\qquad r\ge 0.\label{Kb1}
\end{equation}

If $r\neq 0$, we can use the following formula in
the analysis of (\ref{2I}):
\[
\eqalign{
\int_1^\infty \cos cx\left({\sqrt{1-x^{-2}}\over 1+(r^2-1)x^{-2}}-1
\right)dx=\\
{\sin c \over c}-{1\over
  c}\int_1^\infty\sin(cx)\frac{(x^2-1)(r^2-1)+r^2x^2}
{\sqrt{x^2-1}(x^2+r^2-1)^2}dx}
\]
(integration by parts). In this case we get 
expression (\ref{K2}) for the kernel. For $r=1$ it takes the form
\begin{equation}
K_{BS}(1,z)=
{1\over\pi z}\int_0^\infty\sin(z\cosh
t)\cosh^{-2}t dt.\label{Kb2}
\end{equation}
Recalling the integral representation of Bessel functions, we see that
\[
\frac{d^2}{dz^2}\left(zK_{BS}(1,z)\right)=-{J_0(z)\over 2}, \qquad z>0.
\]
Multiply this relation by $z$ and integrate from zero to $z$. We
integrate the l.h.s. by parts, and for the r.h.s. use the result 
$\int_0^z J_0(z)zdz=zJ_1(z)$. This leads to the following 
equation:
\[
\frac{d}{dz}K_{BS}(1,z)=-{J_1(z)\over 2z}, \qquad z>0.
\]
Its solution satisfying the condition $K_{BS}(\infty)=0$ 
(see (\ref{Kb2})) gives:
\[
K_{BS}(1,z)=\int_z^\infty{J_1(t)\over 2t}dt, \qquad z>0.
\]
Since by (\ref{Kb2}) $K_{BS}(1,z)=K_{BS}(1,-z)$ this completes the proof.
$\Box$

\noindent
{\it Second proof for $r=1$.} 
In this case, let us write the weight function (\ref{w2}) in the form
\[
w(x)(1-\gamma^2x^2)=\sqrt{1-x^2}.
\]
We see that the polynomials $Q_k(x)$ are the Chebyshev polynomials of
the second kind. We have
\begin{equation}
Q_k(x)=\frac{\sin(k+1)\psi}{2^k\sin\psi},\quad x=\cos\psi,\quad
h'_k={\pi\over 2^{2k+1}},\quad k=0,1,\dots\label{h'}
\end{equation}
In the variable $\psi$, the quantity $t_n$ defined by (\ref{phiQ})
has the form
\begin{equation}
t_n=\frac{\gamma}{2^n\pi}\int_0^\pi\frac{
\cos(n\;{\rm arccos}\;(\gamma\cos\psi))}{1-\gamma^2\cos^2\psi}
\sin(n+1)\psi \sin\psi d\psi.\label{tas}
\end{equation}
The asymptotics of this integral as $\alpha=2s/n$ ($\gamma=\cos(\alpha/2)$)
and $n\to\infty$
are analyzed similarly to those of (\ref{Ik}) in the proof of Theorem
3.2. The main contribution to the integral comes from 
the small neighbourhoods of the points zero and $\pi$. The result is
\[
t_n=2^{-n}\left(\frac{2}{\pi}\int_0^\infty
\frac{\cos\sqrt{x^2+s^2}}{x^2+s^2}x\sin x dx+{1\over2}\right)(1+o(1)).
\]
Representing the product of cosine and sine functions as a sum of two
sines and changing the variables $y=x\pm\sqrt{x^2+s^2}$, we finally
obtain
\[
t_n=2^{-n}\left({1\over\pi}\int_0^\infty{dy\over y}{y^2-s^2 \over y^2+s^2}
\sin y+{1\over2}\right)(1+o(1))=2^{-n}e^{-s}(1+o(1)).
\]
Substituting this and $h'_k$ (\ref{h'}) into (\ref{D2}), we obtain
(\ref{Db}) for $r=1$. The rest of the argument is the same as in the first
proof. $\Box$

\section{Toeplitz determinants for $\sin^{\pm1}(\theta/2)$ on an arc}
\noindent
{\bf Lemma 4.1} {\it Let $f_1(\theta)=\sin\theta/2$, 
$f_2(\theta)=\sin^{-1}\theta/2$
if $\alpha\le\theta\le 2\pi-\alpha$, and $f_1(\theta)=f_2(\theta)=0$ 
otherwise. Then for $\alpha=2s/n$, large $n$, we have 
\begin{equation}
\eqalign{
D_{n-1}(f_1)=2^{-n}n^{1/4}
\left(\sqrt{\pi}2^{1/12}e^{3\zeta'(-1)}e^{-s^2/2}J_0(is)+o(1)\right),\\
D_{n-1}(f_2)=2^{n-1}n^{1/4}
\left(\sqrt{\pi}2^{1/12}e^{3\zeta'(-1)}e^{-s^2/2}F(s)
+o(1)\right),\\
F(s)={1\over\pi}\int_0^\infty\frac{\cos\sqrt{x+s^2}}
{x+s^2}J_0(\sqrt{x})dx,}
\end{equation}
where $\zeta'(x)$ is the derivative of Riemann's zeta function, and
$J_0(x)$ is the 0th order Bessel function.

\noindent
Proof.} Combining Lemmas 2.2. and 3.1 (\ref{D}), we obtain:
\begin{equation}
D_{n-1}(f_1)=\frac{2^{n^2-n}\gamma^{n^2+n}}{\pi^n}A_nL_n(1/\gamma),
\qquad A_n=\prod_{j=0}^{n-1}\frac{2^{2j}}{(j+1/2){2j\choose j}^2}.
\label{Dprim}
\end{equation}
Here $\gamma(n)=\cos(s/n)$. For the Legendre polynomial
we use the well-known Hilb asymptotics which is uniform in $s/n$:
\[
 L_n(1/\gamma)={2^n\over {2n\choose n}}(J_0(is)+o(1))=
 {\sqrt{\pi n}\over 2^n}(J_0(is)+o(1)).
\]
The asymptotics of $A_n$ are easy to compute (see \cite{WidomArc}):
\[
A_n=2^{1/12}e^{3\zeta'(-1)}n^{-1/4}(2\pi)^n2^{-n^2}(1+o(1)).
\]
Substituting these results into (\ref{Dprim}) gives the first
asymptotics in the lemma. 

The determinant $D_n(f_2)$ is obtained by using the
second part of Lemma 2.2. and (\ref{D2}) of Lemma 3.1.:
\begin{equation}
D_n(f_2)=(2\gamma)^{n^2+n}{t_n\over \pi^n}A_n,\label{D2prim}
\end{equation}  
and to get its asymptotics, it only remains to calculate $t_n$, 
$n\to\infty$. 
Exactly as (\ref{tas}) we have
\[
t_n=\frac{\gamma}{\pi}\int_0^\pi\frac{
\cos(n\;{\rm arccos}\;(\gamma\cos\psi))}{1-\gamma^2\cos^2\psi}
L_n(\cos(\psi))\sin\psi d\psi.
\]
Now we use Hilb's asymptotics for $L_n(\cos(\psi))$ and estimate the
integral in the same way as (\ref{tas}). We get
\[
t_n={\sqrt{\pi n}\over 2^n}\left({1\over\pi}\int_0^\infty
\frac{\cos\sqrt{x+s^2}}{x+s^2}J_0(\sqrt{x})dx+o(1)\right),
\]
which, in view of (\ref{D2prim}), completes the proof.$\Box$

\section{Acknowledgements}
I am grateful to L. Golinskii and B. Simon for some useful references
regarding polynomials on an arc. This work was supported by the Sfb 288.

\end{document}